\documentstyle[12pt]{article} 
\begin{document}
\newcommand{\singlespace}{
    \renewcommand{\baselinestretch}{1}
\large\normalsize}
\newcommand{\doublespace}{
   \renewcommand{\baselinestretch}{1.2}
   \large\normalsize}
\renewcommand{\theequation}{\thesection.\arabic{equation}}

\input amssym.def
\input amssym
\def \ten#1{_{{}_{\scriptstyle#1}}}
\def \Z{\Bbb Z}
\def \C{\Bbb C}
\def \R{\Bbb R}
\def \Q{\Bbb Q}
\def \N{\Bbb N}
\def \l{\lambda}
\def \L{\Lambda}
\def \M{\Bbb M}
\def \mm{V^{\natural}}
\def \V{V^{\natural}}
\def \wt{{\rm wt}}
\def \tr{{\rm tr}}
\def \Res{{\rm Res}}
\def \End{{\rm End}}
\def \Aut{{\rm Aut}}
\def \mod{{\rm mod}}
\def \Hom{{\rm Hom}}
\def \im{{\rm im}}
\def \cf{{\cal F}}
\def \cs{{\cal S}}
\def \cm{{\cal M}}
\def \D{{\Delta}}
\def \bc{{\bf 1}}
\def \ci{{\cal I}}
\def \<{\langle} 
\def \>{\rangle} 
\def \w{\omega}
\def \o{\omega}
\def \t{\tau }
\def \char{{\rm char}}
\def \a{\alpha }
\def \b{\beta}
\def \e{\epsilon }
\def \la{\lambda }
\def \om{\omega }
\def \O{\Omega}
\def \qed{\mbox{ $\square$}}
\def \pf{\noindent {\bf Proof: \,}}
\def \voa{vertex operator algebra\ }
\def \voas{vertex operator algebras\ }
\def \p{\partial}

\singlespace
%\doublespace
\newtheorem{thmm}{Theorem}
\newtheorem{pp}[thmm]{Proposition}
\newtheorem{th}{Theorem}[section]
\newtheorem{prop}[th]{Proposition}
\newtheorem{coro}[th]{Corollary}
\newtheorem{lem}[th]{Lemma}
\newtheorem{rem}[th]{Remark}
\newtheorem{de}[th]{Definition}

\begin{center}
{\Large {\bf Monstrous Moonshine of higher weight}} \\
\vspace{0.5cm}

Chongying Dong\footnote{Supported by NSF grant 
DMS-9700923 and a research grant from the Committee on Research, UC Santa Cruz.} and Geoffrey Mason\footnote{Supported by NSF grant DMS-9700909 and a research grant from the Committee on Research, UC Santa Cruz.}\\
Department of Mathematics, University of
California, Santa Cruz, CA 95064 \\
\end{center}
\hspace{1.5 cm}

\begin{abstract} We determine the space of 1-point correlation 
functions associated with the Moonshine module: they are precisely those
modular forms of non-negative integral weight which are holomorphic in
the upper half plane, have a pole of order at most 1 at infinity, and
whose Fourier expansion has constant $0.$ There are
Monster-equivariant analogues in which one naturally associates to
each element $g$ in the Monster a modular form of fixed weight $k,$
the case $k=0$ corresponding to the original ``Moonshine'' of Conway
and Norton.
\end{abstract}

\section{Introduction}

Suppose that $V$ is a vertex operator algebra. One of the basic problems is 
that of determining the so-called $n$-point correlation functions 
associated to $V.$ There is a recursive procedure whereby $n$-point functions
determine $n+1$-point functions (see [Z], for example), so that 
understanding 1-point functions become important. In this paper
we will study the 1-point functions (on the torus) associated with
the Moonshine module, which is of interest not only as an example of the
general problem but because of connections with 
the Monster simple group $\M.$

First we recall the definition of a 1-point function. Let the decomposition
of $V$ into homogeneous spaces be given by
\begin{equation}\label{1.1}
V=\bigoplus_{n\geq n_0}V_n.
\end{equation}
Each $v\in V$ is associated to a vertex operator 
\begin{equation}\label{1.2}
Y(v,z)=\sum_{n\in\Z}v(n)z^{-n-1}
\end{equation}
with $v(n)\in \End V.$ If $v$ is homogeneous of weight $k,$ that is
$v\in V_k,$ we write $\wt v=k.$ The {\em zero mode} of $v$ is defined
for homogeneous $v$ to be the component operator 
\begin{equation}\label{1.3}
o(v)=v(\wt v-1)
\end{equation}
and one knows that $o(v)$ induces an endomorphism of each homogeneous 
space. That is,
\begin{equation}\label{1.4}
o(v): V_n\to V_n.
\end{equation}
The 1-{\em point function determined by} $v$ is then essentially the graded trace of $o(v)$ on $V.$ More precisely,
if $V$ has central charge $c$ we define the 1-point function 
(on the torus) via
\begin{equation}\label{1.5}
Z(v,q)=Z(v,\tau)=\tr|_Vo(v)q^{L(0)-c/24}=q^{-c/24}\sum_{n\geq n_0}(\tr|_{V_n}o(v))q^n.
\end{equation}
Here, $L(0)$ is the usual degree operator and $q$ may be taken either
as an indeterminate or, less formally, to be $e^{2\pi i \tau}$ with
$\tau$ in the upper half plane $\frak h.$ If $g$ is an automorphism
of $V$ we define
\begin{equation}\label{1.6}
Z(v,g,q)=Z(v,g,\tau)=q^{-c/24}\sum_{n\geq n_0}(\tr|_{V_n}o(v)g)q^n.
\end{equation}
These functions can be extended linearly to all $v\in V$ by defining
$Z(v,g,q)=\sum_{i}Z(v_i,g,q)$ if $v=\sum_{i}v_i$ is the decomposition
of $v$ into homogeneous components. In this way we obtain
the space of 1-point functions associated to $V,$ namely the functions 
$Z(v,q)$ for $v\in V.$ 

In order to state our results efficiently we need some notation concerning
modular forms. We denote ${\cal F}$ the $\C$-linear space spanned
by those (meromorphic) modular forms $f(\tau)$ of level 1 and integral weight 
$k\geq 0$ which satisfy

\bigskip
(i) $f(\tau)$ is holomorphic in $\frak h.$

(ii) $f(\tau)$ has Fourier expansion of the form
\begin{equation}\label{1.7}
f(\tau)=\sum_{n=-1}^{\infty}a_nq^n,\ \  a_0=0.
\end{equation}

\bigskip

Thus $f(\tau)$ has a pole of order at most 1 at infinity
and constant 0. Let ${\cal M}$ be the space of holomorphic
modular forms of level 1 and ${\cal S}$ the space of cusp-forms
of level 1. Thus we have ${\cal S}={\cal F}\cap {\cal M}.$

Among the elements of $\cm$ are the Eisenstain series $E_k(\t)$ for 
even $k\geq 4.$ We normalize them as in [DLM], namely
\begin{equation}\label{1.8}
E_k(\t)=\frac{-B_k}{k!}+\frac{2}{(k-1)!}\sum_{n=1}^{\infty}\sigma_{k-1}(n)q^n
\end{equation}
with $B_k$ the $k$th Bernoulli number defined by 
\begin{equation}\label{1.9}
\frac{t}{e^t-1}=\sum_{k=0}^{\infty}B_k\frac{t^k}{k!}.
\end{equation} 
If $\cm_k$ is the space of forms $f(\t)\in\cm$ of weight $k$ then
there is a differential operator $\partial: \cm_k\to \cm_{k+2}$ defined
via
\begin{equation}\label{1.10}
\partial=\partial_k: f(\tau)\mapsto \frac{1}{2\pi i}\frac{d}{d\tau}f(\tau)
+kE_2(\tau)f(\tau).
\end{equation}
Here, $E_2(\tau)$ is again defined by (\ref{1.8}), though $E_2$ is
not a modular form.

By a $\p$-{\em ideal} we mean an ideal $\ci$ in the commutative algebra 
$\cm$ which also satisfies $\p(\ci)\subset \ci.$
\begin{thmm}\label{thmm1}
Let $V^{\natural}$ be the Moonshine module. The space of 1-point
functions associated to $\mm$ is precisely the linear space $\cf$ 
defined above.
\end{thmm}

As we will explain in due course, it is a consequence of results 
in [Z] (see also [DLM]) that all 1-point functions associated to vectors 
$v\in \mm$ lie in $\cf.$ The new result here is therefore an existence
result: for each $f(\t)\in\cf$ there is a $v\in \mm$ such that $Z(v,\t)
=f(\t).$

Recall next that $\mm$ is a direct sum of irreducible highest weight
modules $M(c,k)$ for the Virasoro algebra $Vir.$ Here, $c=24$ and
for $k>0,$ $M(c,k)$ is the Verma module generated a highest weight 
vector $v\in \mm_k.$ Thus $L(n)v=0$ for
all $n>0$ where $L(n)$ are the usual generators for $Vir,$ and 
$L(0)v=kv.$ 

The proof of Theorem \ref{thmm1} is facilitated by the next 
result.
\begin{pp}\label{p2} Let $v\in\mm_k$ be a highest weight 
vector of positive weight $k.$ Then the following 
hold:

(a) $Z(v,\tau)$ is a cusp-form  of weight $k,$ possibly 0.

(b) The space of 1-point functions consisting of all
$Z(w,\t)$ for $w$ in the highest weight module
for $Vir$ generated by $v$ is the $\p$-ideal
generated by $Z(v,\t).$
\end{pp}

While Proposition \ref{p2} actually holds for a wide class of vertex 
operator algebras, our final result is more 
closely tied to the structure of $\mm.$ It gives
us a large set of highest weight vectors (for the Virasoro algebra)
to which we can usefully apply the preceding proposition.

First recall that to each $\l$ in the Leech Lattice $\L$ there is
a corresponding element $e^{\l}$ in the group algebra $\C[\L]$ and
an element, also denoted $e^{\l},$ in the vertex operator algebra
$V_{\L}$ associated to $\L.$ See [B1], [FLM] and Section 4
below for more details. The relation of $V_{\L}$ to $\mm$
shows that $e^{\l}+e^{-\l}$ can be considered as an element 
of both vertex operator algebras.

\begin{thmm}\label{thmm3}
Let $v(\l)=e^{\l}+e^{-\l}$ be as above and considered as an element of $\mm.$
Then $v(\l)$ is a highest weight vector of weight $k=\frac{\<\l,\l\>}{2}$
and if $0\ne \l\in 2\L$ then 
\begin{equation}\label{1.11}
Z(v(\l),\t)=\eta(\t)^{12}\left\{
\left(\frac{\Theta_1(\tau)}{2}\right)^{\<\l,\l\>-12}
+\left(\frac{\Theta_2(\tau)}{2}\right)^{\<\l,\l\>-12}
-\left(\frac{\Theta_3(\tau)}{2}\right)^{\<\l,\l\>-12}\right\}.
\end{equation}
\end{thmm}

In (\ref{1.11}), $\eta(\t)$ is the Dedekind eta function and $\Theta_1,
\Theta_2, \Theta_3$ are the usual Jacobi theta functions
(see, for example [C], p. 69).

If $\L_n=\{\l\in\L|\<\l,\l\>=2n\}$ then $\L_2=0,$ so if $0\ne \l\in
2\L$ then $\frac{1}{2}\<\l,\l\>=4m$ with $m\geq 2.$ If $m=2$ then
$Z(v(\l),\t)$ is a cusp form of level 1 and weight 8 by Proposition
\ref{p2}, and hence must be 0. Then (\ref{1.11}) reduces to the identity
$\Theta_1(\t)^4+\Theta_2(\t)^4-\Theta_3(\t)^4=0,$ 
which is well-known in the theory of elliptic functions (loc. cit.). 
If $m\geq 3$ then one can check that
$Z(v(\l),\t)\ne 0$ (e.g. by looking at the coefficient of $q$ in the
Fourier expansion), so $Z(v(\l),\t)$ is a non-zero cusp form of level
1 and weight $4m=12, 16, 20,\cdots.$ One knows (see, for example [S])
 that the cusp forms of level 1 and weights $12, 16, 20$ are unique
up to scalar (as are those of weight $18, 22$ and $26$) and given by
$\Delta(\tau),$ $\D(\t)E_4(\t), \D(\t)E_8(\t)$ respectively, where
$\D(\t)=\eta(\t)^{24}$ is the discriminant. Once we know that
$\Delta(\t)$ can be realized as a 1-point function $Z(v,\t)$ for some
highest weight vector $v,$ the fact that $\cs=\cm\D(\t)$ (loc. cit.)
together with Proposition \ref{p2} then shows that every $f(\t)\in\cs$
can be so realized. This in turn reduces the proof of Theorem
\ref{thmm1} to dealing with forms which have a pole at infinity.

Our discussion so far has not taken into account the automorphisms $g$
of $\mm$ (i.e., elements of the Monster). There are some general 
results, which follow from [DLM], which imply that if $v\in \mm$ is homogeneous
of weight $k$ with respect to a certain operator $L[0],$ then 
$Z(v,g,\tau)$ is a modular form of weight $k$ for each $g\in {\Bbb M}.$
Moreover the level is the same as that for the McKay-Thompson series 
$Z({\bf 1},g,\t)$ described in [CN]. We describe the precise subgroup of
$SL(2,\Z)$ which fixes $Z(v,g,\t)$ in Theorem \ref{t6.1}.

Group theorists may be disappointed to learn that if we fix $v$ so
that all $Z(v,g,\t)$ are modular forms of weight $k$ then in general
the Fourier coefficients of the forms (for varying $g$) do {\em not}
define characters, or even generalized characters. This is so even if
$Z(v,1,\t)$ has integer coefficients. This does not mean, however, that these
higher weight McKay-Thompson series are of no arithmetic 
interest. If we combine our results with some
calculations of Harada and Lang [HL], for example, we find that for each 
of the weights $k=12, 16, 20$ there is a {\em unique} vector
$v$ in the Moonshine module $\mm$ with the following properties:

\bigskip

(a) $v$ is a highest weight vector for $Vir$ which lies in $\mm_k$ 
and is Monster-invariant.

(b) The 1-point function $Z(v,\t)=q+\cdots$ is the unique normalized cusp form of level 1 and weight $k.$

\bigskip

\noindent Such a $v$ may be obtained by averaging the vector $v(\l)$ of 
Theorem \ref{thmm3} over the Monster ($\l\in 2\L_m, m=3,4$ or
$5$). The unicity of such
$v$ makes them entirely analogous to the vacuum vector $\bc,$ and it 
is likely that the trace functions $Z(v,g,\t)$ are of
particular interest in these cases.

We can understand the representation-theoretic meaning of the functions 
$Z(v,g,\t)$ as follows: since $v$ is  Monster-invariant then each $g$
commutes with the zero mode $o(v)$ and its semi-simple part $o(v)_s$
with regard to its action on the homogeneous space $\mm_n.$ Thus if $o(v)_s$
has distinct eigenvalues $\l_1,...,\l_t$ on $\mm_n,$ the corresponding
 eigenspaces $\mm_{n,1},...,\mm_{n,t}$ are Monster-modules and the 
$(n-1)th$ Fourier coefficient of $Z(v,g,\t)$ is equal to
$\sum_{i=1}^t\l_i\tr_{\mm_{n,i}}g.$

We complete our discussion with two conjectures: (A) For each cusp
form $f(\t)\in\cs$ of weight $k$ there is a (Monster-invariant)
highest weight vector $v\in \mm_k$ with $Z(v,\t)=f(\t);$ (B) If
$Z(v,\t)$ is a cusp form then so is $Z(v,g,\t)$ for each Monster
element $g.$

The paper is organized as follows: In Section 2 we review the required results
from the theory of vertex operator algebra and prove Proposition
\ref{p2}. In Section 3 we reduce the proof of Theorem \ref{thmm1} to
that of Theorem \ref{thmm3}, which is proved by lengthy calculation in
Section 4. In Section 5 we give an equivariant version of formula
(\ref{1.11}), that is, we calculate $Z(v(\l),g,\t)$ for various (but
not all!) elements $g\in \M,$ and in Section 6 we describe the invariance 
group of $Z(v,g,\t)$ in $SL(2,\Z).$

We thank Chris Cummins for useful comments on a prior version of this 
paper.

\section{Proof of Proposition 2}
\setcounter{equation}{0}

We start by recalling some results from [Z] and [DLM].  If $V$ is
a vertex operator algebra as in (\ref{1.1}) then there is a second 
VOA structure $(V,Y[])$ defined on $V$ with vertex operator
$Y[v,z]$. The two VOAs are related by a change of variables and 
have the same vacuum vector ${\bf 1}$ and central charge $c.$ The conformal
vectors are distinct, however, and we denote the
standard Virasoro generators for the second VOA by $L[n].$ The relation
between the $L(n)$ and $L[n]$ (cf. [Z]) shows that both Virasoro algebras
have the {\em same} highest weight vectors $v.$ 

A most important identity for us is the following (cf. [Z] and [DLM],
equation (5.8)): if $w\in V$ then
\begin{equation}\label{2.1}
Z(L[-2]w,\t)=\p Z(w,\t)+\sum_{l=2}^{\infty}E_{2l}(\tau)Z(L[2l-2]w,\t)
\end{equation}
where we are using the notation of Section 1. We should emphasize that 
it is a consequence of the main results of [Z] and [DLM] that
if $v$ is homogeneous of weight $k$ with  respect to $L[0]$, where we
are taking $V=\V$ to be the Moonshine module, then the trace function 
$Z(v,\tau)$ is indeed a meromorphic modular form of level 1 which lies in the space $\cf$ defined in (\ref{1.7}).

It is also shown in [Z] (cf. [DLM], equation (5.1)) that the
following holds:
\begin{equation}\label{2.2}
Z(L[-1]w,\t)=0\ \ \ {\rm for \ all}\ w\in V.
\end{equation}

We turn to the proof of Proposition \ref{p2}, beginning with part (a), which 
is elementary. Namely, from the creation axiom 
$$\lim_{z\to 0}Y(v,z)\bc=v$$
we get $v(n)\bc=0$ if $n\geq 0.$ So if $v\in V_k$ with $k>0$ then 
$o(v)\bc=0,$ in which case we see that  
$$Z(v,\t)=q^{-1}\sum_{n=2}^{\infty}\tr|_{V_n}o(v)q^n$$
is a modular form of level 1, holomorphic in $\frak h$ with a zero of order 
at least 1 at $\infty.$ So 
indeed $Z(v,\t)$ is a cusp-form, as asserted in Proposition \ref{p2} (a).

We turn to the proof of (b) of Proposition \ref{p2}, which is established 
by a systematic use of equations (\ref{2.1}) and (\ref{2.2}). Let $v\in V_k$ be
a highest weight vector. By a {\em descendant} of $v$ we will mean 
a vector of the form $L[n_1]\cdots L[n_t]v$
with each $n_i\leq 0,$ or any linear combination of such vectors; we write 
$v\to w$ if $w$ is a descendant of $v.$  

Let $I=\<Z(w,\tau)|v\to w\>$ be the linear span of the indicated forms
and let $J$ be the $\p$-ideal generated by 
$Z(v,\t).$ We must prove that $I=J.$

First we show that $I\subset J.$ We do this by proving by induction
on wt[$w$] (the weight of $w$, homogeneous with respect to the second Virasoro
algebra) that $Z(w,\t)\in J.$ For this we may take $w$ in the 
form $w=L[n_1]\cdots L[n_t]v$ with each $n_i=-1$ or $-2.$ If $n_1=-1$ then 
$Z(w,\t)=0$ by (\ref{2.2}), so we may take $n_1=-2.$ So 
$w=L[-2]x$ where $x=L[n_2]\cdots L[n_t]v$ has weight equal to wt[$w$]-2.

By (\ref{2.1}) we have 
\begin{equation}\label{2.3}
Z(w,\t)=\p Z(x,\t)+\sum_{l=2}^{\infty}E_{2l}(\tau)Z(L[2l-2]x,\t).
\end{equation}
Since $v\to x$ and $v\to L[2k-2]x,$ induction tells us that $Z(x,\t)$
and $Z(L[2l-2]x,\t)$ both lie in $J,$ whence so does of r.h.s. of (\ref{2.3})
since $J$ is a $\p$-ideal. So indeed $Z(w,\t)$ lies in $J.$

Next we show that $I$ is also a $\p$-ideal. Since $Z(v,\t)$ is in $I$
it follows from this that $J\subset I$ and hence that $I=J,$ as required.

Let $r\geq 1$ with $v\to w$ and consider the vector
$x=L[-2]L[-1]^{2r}w.$ If $2l-2<2r$ then $L[2l-2]L[-1]^{2r}w$ can be
written as a linear combination of vectors of the shape $L[-1]u$ for
some $u.$ Thus (\ref{2.2}) tells us that
$Z(L[2l-2]L[-1]^{2r}w,\tau)=0$ if $2l-2<2r.$ Now by (\ref{2.1}) we get
\begin{equation}\label{2.4}
Z(x,\t)=\sum_{l=r+1}^{\infty}E_{2l}(\tau)Z(L[2l-2]L[-1]^{2r}w,\t).
\end{equation}
Assuming that $w$ is homogeneous with respect to the second Virasoro
algebra, it follows in the same way that $Z(L[2r]L[-1]^{2r}w,\t)$ is a
non-zero multiple of $Z(w,\t).$ If $l>r+1$ then
$L[2l-2]L[-1]^{2r}w$ has weight less than that of $w,$ while if also 
$v=w$ then  $L[2l-2]L[-1]^{2r}v=0.$ Thus (\ref{2.4}) now reads
\begin{equation}\label{2.5}
Z(x,\t)=\alpha E_{2r+2}(\t)Z(w,\tau)+\sum_{l=r+2}^{\infty}E_{2l}(\tau)Z(u_l,\t)
\end{equation}
where $v\to u_l,$ wt[$u_l$]$<$wt[$w$] and $\a$ is a non-zero
scalar. From (\ref{2.5}) and what we have said it follows by induction
on wt[$w$] that $E_{2r+2}(\t)Z(w,\t)$ lies in $I$ whenever $r\geq 1.$
Since the forms $E_{2r+2}(\t)$ generate the space $\cm$ of modular
forms (in fact $E_4(\t)$ and $E_6(\t)$ suffice), it follows that $I$
is an ideal in $\cm.$ But then (\ref{2.1}) shows that $\p Z(w,\t)$
lies in $I$ whenever $v\to w,$ so $I$ is a $\p$-ideal. This completes
the proof of Proposition \ref{p2} (b).

\section{Trace functions with a  pole}
\setcounter{equation}{0}

In this section we prove
\begin{prop}\label{p3.1}
Let $k$ be a non-negative integer. Then the trace function
$Z(L[-2]^k\bc,\tau)$ is non-zero, and more precisely has a
$q$-expansion of form $\epsilon q^{-1}+\cdots $ where $(-1)^k\epsilon>0.$
\end{prop}

Set $w=L[-2]^k\bc.$ Note that the truth of the proposition shows that
$Z(w,\tau)$ is a form of level 1 and weight $2k$ which is non-zero
with a pole at $\infty.$ If we have two such trace functions of the
same weight and the same residue at $\infty$ then they differ by a
cusp-form. So together with Proposition \ref{p2} , this reduces the
proof of Theorem \ref{thmm1} to showing that $\Delta(\tau),$ say, can be realized
as a trace function. As we have pointed out in Section 1, this is
implicit in the statement of Theorem \ref{thmm3}.

We turn to the proof
of Proposition
\ref{p3.1}, using induction on $k.$ The case $k=0$ is obvious. Set 
$x=L[-2]^{k-1}\bc, $ so that $w=L[-2]x.$ By  (\ref{2.1}) and (\ref{1.10})
we get
\begin{equation}\label{3.1}
Z(w,\t)=q\frac{d}{dq}Z(x,\t)+\sum_{l=1}^{\infty}E_{2l}(\tau)Z(L[2l-2]x,\t).
\end{equation}

Now by another induction argument using the Virasoro relations, we
easily find that if $l\geq 1$ then there is an identity of
the form 
\begin{equation}\label{3.2}
L[2l-2]x=n_lL[-2]^{k-l}\bc
\end{equation}
where $n_l$ is positive and the right-side is interpreted as $0$ if $l>k.$

From (\ref{1.8}), the $q$-expansion of $E_{2l}(\tau)$ begins 
$$-\frac{B_{2l}}{(2l)!}+\cdots$$
and it is easily seen from (\ref{1.9}) that we have
\begin{equation}\label{3.3}
(-1)^{l+1}B_{2l}>0.
\end{equation}

By induction we have $Z(L[-2]^r\bc,\t)=\epsilon(r)q^{-1}+\cdots $ with
$(-1)^r\epsilon(r)>0$ for $0\leq r<k.$ It follows that the coefficient
of $q^{-1}$ on the r.h.s. of (\ref{3.1}) is equal to
\begin{eqnarray*}
& &\ \ \ \ -\e(k-1)-\sum_{l=1}^k\frac{B_{2l}}{(2l)!}n_l\e(k-l)\\
& &=(-1)^k\left\{(-1)^{k-1}\e(k-1)+\sum_{l=1}^k(-1)^{l+1}\frac{B_{2l}}{(2l)!}n_l(-1)^{k-l}\e(k-l)\right\}.
\end{eqnarray*}
From what we have said, the sum of the terms in the braces is positive,
so Proposition \ref{3.1} is proved. \qed

\section{Proof of Theorem \ref{thmm3}}
\setcounter{equation}{0}

We have reduced the proof of Theorem \ref{thmm1} to that
of Theorem \ref{thmm3}, which we carry out in this section. 

We first take over {\it en bloc}  the notation of [FLM] with regard to  the
lattice VOA $V_{\L}$ and associated vertex operators, where $\L$ is the
Leech lattice. In particular, $\frak h=\C\otimes_{\Z}\L;$ $\hat{\frak
h}_{\Bbb Z}$ is the corresponding Heisenberg algebra; $M(1)$ is the
associated irreducible induced module for $\hat{\frak h}_{\Bbb Z}$
such that the canonical central element of $\hat{\frak h}_{\Bbb Z}$
acts as 1; $V_{\L}=M(1)\otimes {\Bbb C}[\L];$
$$Y(e^\alpha,z)=E^-(-\a,z)E^+(-\a,z)e_{\alpha}z^{\a}$$ 
is the vertex
operator associated to $\alpha\in \L$ where
$$E^{\pm}(\alpha,z)=\exp\left(\sum_{n\in \N}\frac{\a(\pm n)}{\pm
n}z^{\mp n}\right)$$ for $\a\in\frak h$ and $e_{\a}$ acts on $\C[\L]$
by $$e_{\alpha}: e^{\beta}\mapsto \e(\alpha,\beta)e^{\a+\b}$$
where $\e(\cdot,\cdot)$ is a bilinear 2-cocycle of $\L$ with values in
$\{\pm 1\}$; $t$ is the automorphism of $V_{\L}$ of order 2 induced
from the $-1$ isometry of $\L$ such that $t e^{\alpha}=e^{-\a};$ 
$t$ acts on $M(1)$ by
$t(\b_1(-n_1)\cdots \b_k(-n_k))=(-1)^k\b_1(-n_1)\cdots \b_k(-n_k)$
for $\b_i\in\frak h$ and $n_i>0.$ 

For a $t$-stable subspace $W$ of $V_{\L}$ we define $W^{\pm}$ to be the
eigenspaces of $t$ with eigenvalues $\pm 1.$ 
We start by considering the action of $Y(e^{\a}+e^{-\a},z)$ on $V^+_{\L}.$
Thus $V_{\L}^+$ is spanned by elements of the form
\begin{equation}\label{4.1}
v\otimes e^{\b}+tv\otimes e^{-\b}
\end{equation} 
and we have 
\begin{eqnarray}
& &\ \ \  Y(v(\a),z)(v\otimes e^{\b}+tv\otimes e^{-\b})\nonumber\\
& &=z^{\<\a,\b\>}E^-(-\a,z)E^+(-\a,z)v\otimes \e(\a,\b)e^{\a+\b}\nonumber\\
& &\ \ \ \ \ \ \ +
z^{-\<\a,\b\>}E^-(-\a,z)E^+(-\a,z)tv\otimes \e(\a,-\b)e^{\a-\b}\nonumber\\
& &+z^{-\<\a,\b\>}E^-(\a,z)E^+(\a,z)v\otimes \e(-\a,\b)e^{-\a+\b}\nonumber\\
& &\ \ \ \ \ \ \ +
z^{\<\a,\b\>}E^-(\a,z)E^+(\a,z)tv\otimes \e(-\a,-\b)e^{-\a-\b}\label{4.2}
\end{eqnarray}
From this we see that non-zero contributions to the trace on $V_{\L}^+$
can arise only when $\alpha\in 2\L,$ and more precisely when $\alpha=\pm
2\beta$ in (\ref{4.2}).  

For $\b\in \L$ we set
$$V(\b)=M(1)\otimes(\C e^{\b}+\C e^{-\b})$$
which is $t$-stable. So  the trace of $o(v(\l))$ on $V_{\L}^+$ is 
equal to the trace of $o(v(\l))$ on $V(\a)^+$ where 
$\l=2\a\in \L,$ which we now assume. Clearly $v(\l)$ is a highest weight
vector with weight $\frac{\<\l,\l\>}{2}.$ 

Note that $\e(\pm\l,\pm\a)=1.$ It follows from (\ref{4.2}) that
 only expressions of the form 
\begin{equation}\label{4.3}
z^{-2\<\a,\a\>}\left(E^-(-2\a,z)E^+(-2\a,z)tv\otimes e^{\a}
+E^-(2\a,z)E^+(2\a,z)v\otimes e^{-\a}\right)
\end{equation}
contribute to the trace. Thus we are essentially reduced 
to computing the trace of the 
degree zero operators of $E^-(-2\a,z)E^+(-2\a,z)$ and $E^-(2\a,z)E^+(2\a,z)$ on
$M(1).$

Let $A=\C\l$ and $\frak h=A\bot B$ be an orthogonal direct sum. Then 
$M(1)=S(\hat \frak h^-)=S(\hat A^-)\otimes S(\hat B^-).$ Let $x$ be 
a formal variable and define $x^N\in (\End M(1))[x]$ such that
$x^N(\a_1(-n_1)\cdots \a_k(-n_k))=x^k\a_1(-n_1)\cdots \a_k(-n_k)$
for $\a_i\in \frak h$ and $n_i>0.$ Set
$$E^-(\pm\l,z)E^+(\pm\l,z)=\sum_{n\in\Z}E^{\pm}(n)z^{-n}.$$

\begin{lem}\label{l4.1} We have
\begin{equation}\label{4.4}
\tr E^{\pm}(0)q^{L(0)}x^N|_{S(\hat A^-)}=
\frac{\exp\left(\sum_{n>0}\frac{-\<\l,\l\>xq^n}{n(1-xq^n)}\right)}{\prod_{n>0}
(1-xq^n)}.
\end{equation}
\end{lem}

\pf Note that $S(\hat A^-)$ has a  basis 
$$\{\l(-n)^{k_n}\cdots \l(-1)^{k_1}|k_i\geq 0, n\geq 1\}.$$ In order
to compute the trace it suffices to compute the coefficients of
$\l(-n)^{k_n}\cdots \l(-1)^{k_1}$ in $E(0)^{\pm}\l(-n)^{k_n}\cdots
\l(-1)^{k_1}.$ That is, we need to compute the projection 
$$P_{k_1,...,k_n}: E(0)^{\pm}\l(-n)^{k_n}\cdots \l(-1)^{k_1}\to 
\C\l(-n)^{k_n}\cdots \l(-1)^{k_1}.$$
 Recall that
$$[\l(s),\l(t)]=s\<\l,\l\>\delta_{s+t,0}$$
for $s,t\in\Z.$ Then
\begin{eqnarray*}
& &\ \ \ \ \ P_{k_1,...,k_n}E(0)^{\pm}\l(-n)^{k_n}\cdots
\l(-1)^{k_1}\\
& &=\sum_{p_i\leq k_i}(-1)^{p_1+\cdots p_n}\frac{\l(-1)^{p_1}}{p_1!}
\cdots \frac{\l(-n)^{p_n}}{n^{p_n}p_n!}\frac{\l(1)^{p_1}}{p_1!}
\cdots \frac{\l(n)^{p_n}}{n^{p_n}p_n!}\l(-n)^{k_1}\cdots \l(-1)^{k_1}\\
& &=\sum_{p_i\leq k_i}(-1)^{p_1+\cdots p_n}\left(\prod_{i=1}^n\frac{
\<\l,\l\>^{p_i}i^{p_i}k_i(k_i-1)\cdots (k_i-p_i+1)}{(p_i!)^2i^{2p_i}}\right)\l(-n)^{k_1}\cdots \l(-1)^{k_1}\\
& &=\sum_{p_i\leq k_i}\left(\prod_{i=1}^n\frac{{k_i\choose p_i}}{p_i!}
\left(\frac{-\<\l,\l\>}{i}\right)^{p_i}\right)\l(-n)^{k_1}\cdots \l(-1)^{k_1}.
\end{eqnarray*}
Thus 
\begin{eqnarray*}
& &\ \ \ \ \ \tr E^{\pm}(0)q^{L(0)}x^N|_{S(\hat A^-)}\\
& &=\sum_{n\geq 1}\sum_{k_i,p_i\geq 0}\left(\prod_{i=1}^n\frac{{k_i\choose p_i}}{p_i!}\left(\frac{-\<\l,\l\>}{i}\right)^{p_i}\right)q^{k_1+2k_2+\cdots +nk_n}
x^{k_1+\cdots k_n}\\
& &=\prod_{i\geq 1}\left(\sum_{k_i,p_i\geq 0}\frac{{k_i\choose p_i}}{p_i!}\left(\frac{-\<\l,\l\>}{i}\right)^{p_i}q^{ik_i}x^{k_i}\right).
\end{eqnarray*}
Note that if $y$ is a formal variable and $s$ is a nonnegative integer
then
\begin{eqnarray*}
& &\ \ \ \ \ \sum_{m\geq s}{m \choose s}y^m\\
& &=y^s\sum_{m=0}^{\infty}{s+m\choose s}y^m\\
& &=y^s\sum_{m=0}^{\infty}{s+m\choose m}y^m\\
& &=\frac{y^s}{(1-y)^{1+s}}.
\end{eqnarray*}
Then for any $p_i\geq 0$ we have
$$ \sum_{k_i\geq 0}{k_i \choose p_i}(xq^i)^{k_i}=\frac{(xq^i)^{p_i}}{(1-xq^i)^{1+p_i}}.$$
Hence 
\begin{eqnarray*}
& &\ \ \ \ \ \tr E^{\pm}(0)q^{L(0)}x^N|_{S(\hat A^-)}\\
& &=\prod_{i\geq 1}\sum_{p_i\geq 0}\frac{1}{(1-xq^i)}
\frac{1}{p_i!}\left(\frac{-\<\l,\l\>xq^i}{i(1-xq^i)}\right)^{p_i}\\
& &=\prod_{n\geq 1}\frac{1}{(1-xq^n)}\exp\left(\frac{-\<\l,\l\>xq^n}{n(1-xq^n)}\right)\\
& &=\prod_{n=1}^{\infty}\frac{1}{(1-xq^n)}\exp\left(\sum_{n=1}^{\infty}\frac{-\<\l,\l\>xq^n}{n(1-xq^n)}\right),
\end{eqnarray*}
as desired. \qed

\begin{lem}\label{l4.2} We have
\begin{equation}\label{4.5}
\tr E^{\pm}(0)q^{L(0)}x^N|_{M(1)}=
\frac{\exp\left(\sum_{n>0}\frac{-\<\l,\l\>xq^n}{n(1-xq^n)}\right)}{\prod_{n>0}
(1-xq^n)^{24}}.
\end{equation}
\end{lem}

\pf Since $M(1)=S(\hat A^-)\otimes S(\hat B^-)$ and $E^{\pm}(0)$ commute
with $\beta(n)$ for $\beta\in B$ and $n\in \Z,$ we immediately have
$$\tr E^{\pm}(0)q^{L(0)}x^N|_{M(1)}=\tr E^{\pm}(0)q^{L(0)}x^N|_{S(\hat A^-)}
\tr q^{L(0)}x^N|_{S(\hat B^-)}$$
and also 
$$\tr q^{L(0)}x^N|_{S(\hat B^-)}=\frac{1}{\prod_{n>0}
(1-xq^n)^{23}}.$$
The lemma now follows from Lemma \ref{l4.1}.
\qed

\bigskip

Set $f(q,x)=\tr E^{\pm}(0)q^{L(0)}x^N|_{M(1)}.$ Then one can easily see that
\begin{eqnarray}\label{4.6}
& &\tr E^{\pm}(0)q^{L(0)}|_{M(1)^+}=\frac{1}{2}(f(q,1)+f(q,-1))\nonumber\\
& &\tr E^{\pm}(0)q^{L(0)}|_{M(1)^-}=\frac{1}{2}(f(q,1)-f(q,-1)).
\end{eqnarray}

\begin{lem}\label{l4.3} The contribution of $V_{\L}^+$ to $Z(v(\l),\tau)$ is
\begin{equation}\label{4.7}
q^{\frac{1}{8}\<\l,\l\>-1}\prod_{n=1}^{\infty}\frac{(1-q^n)^{24}}{(1-q^{2n})^{24}}\prod_{n=1}^{\infty}\frac{(1-q^{2n})^{2\<\l,\l\>}}{(1-q^{n})^{\<\l,\l\>}}
=\frac{\eta(2\t)^{2\<\l,\l\>-24}}{\eta(\t)^{\<\l,\l\>-24}}.
\end{equation}
\end{lem}

\pf We have already seen that
$$\tr o(v(\l))q^{L(0)}|_{V_{\L}^+}=\tr o(v(\l))q^{L(0)}|_{V(\a)^+}.$$
Clearly, $q^{L(0)}e^{\pm\a}=q^{\frac{1}{8}\<\l,\l\>}e^{\pm \a}.$
From the proof of Lemma \ref{l4.1} we see that $E^{\pm}(0)$ 
have the same eigenvectors and the corresponding eigenvalues are also
the same.  
It follows from (\ref{4.3}), (\ref{4.5}) and (\ref{4.6}) that 
\begin{eqnarray*}
& &\ \ \ \ \ \ \tr o(v(\l))q^{L(0)}|_{V(\a)^+}\\
& &=q^{\frac{1}{8}\<\l,\l\>}(E^{\pm}(0)q^{L(0)}|_{M(1)^+}-\tr E^{\pm}(0)q^{L(0)}|_{M(1)^-})\\
& &=q^{\frac{1}{8}\<\l,\l\>}f(q,-1)\\
& &=q^{\frac{1}{8}\<\l,\l\>}\frac{\exp\left(\sum_{n>0}\frac{\<\l,\l\>q^n}{n(1+q^n)}\right)}{\prod_{n>0}
(1+q^n)^{24}}.
\end{eqnarray*}
Next note that
\begin{eqnarray*}
& &\sum_{n>0}\frac{q^n}{n(1+q^n)}=\sum_{n=1}^{\infty}\frac{q^n}{n}\sum_{i=0}^{\infty}(-1)^iq^{in}\\
& &\ \ \ \ \ =-\sum_{i=1}^{\infty}(-1)^i\sum_{n=1}^{\infty}\frac{q^{in}}{n}\\
& &\ \ \ \ \ = \sum_{n=1}^{\infty}(-1)^n\log(1-q^n).
\end{eqnarray*}
So $\tr o(v(\l))q^{L(0)}|_{V(\a)^+}$ may be written as
\begin{equation}\label{4.8}
q^{\frac{1}{8}\<\l,\l\>}\prod_{n=1}^{\infty}(1+q^n)^{-24}\prod_{n=1}^{\infty}(1-q^n)^{(-1)^n\<\l,\l\>}.
\end{equation}

If now we incorporate the grade-shift of $q^{-c/24}=q^{-1}$, the lemma follows from (\ref{4.8}).
\qed

At this point, recall [C] the Jacobi theta functions $\Theta_i,$ $i=1,2,3,$ 
considered as functions of $\t$ i.e., with the ``other'' variable set equal
to 0: 
\begin{equation}\label{4.9}
\Theta_1(\t)=2q^{\frac{1}{8}}\prod_{n=1}^{\infty}(1-q^n)(1+q^n)^2=
2\frac{\eta(2\t)^2}{\eta(\t)}
\end{equation}
\begin{equation}\label{4.10}
\Theta_2(\t)
=\prod_{n=1}^{\infty}(1-q^n)(1-q^{n-1/2})^2=
\frac{\eta(\t/2)^2}{\eta(\t)}
\end{equation}
\begin{equation}\label{4.11}
\Theta_3(\t)=\prod_{n=1}^{\infty}(1-q^n)(1+q^{n-1/2})^2=
\frac{\eta(\t)^5}{\eta(\t/2)^2\eta(2\t)^2}.
\end{equation}

Combining (\ref{4.7}) and (\ref{4.9}) then yields 
\begin{lem}\label{l4.9} The contribution  
of $V_{\L}^+$ to $Z(v(\l),\tau)$ is equal to 
$$\eta(\t)^{12}\left(\frac{\Theta_1(\t)}{2}\right)^{\<\l,\l\>-12}.$$
\end{lem}

Now let $V_{\L}^T$ be the $t$-twisted $V_{\L}$-module (cf. [FLM]). Then
the moonshine module $V^{\natural}$ is the direct sum of $V_{\L}^+$ and 
$(V^T_{\L})^+$ where again $+$ refer to the fixed points of the action $t$
on $V_{\L}^T.$ The space $V_{\L}^T$ can be described as follows:
$$V_{\L}^T=S(\hat\frak h[-1]^-)\otimes T$$
where $\hat \frak h[-1]=\sum_{n\in\Z}\frak h\otimes t^{n+1/2}\oplus \C c$
is the $-1$-twisted Heisenberg algebra, $ \hat \frak h[-1]^-=\sum_{n>0}
\frak h\otimes t^{-n+1/2}$ and $T$ is the $2^{12}$-dimensional
projective representation for $\L$ such that $2L$ acts on $T$ trivially.
The grading
on $V_{\L}^T$ is the natural one together with an overall shift 
of $q^{3/2}.$ Now $t$ acts on $T$ as multiplication by $-1$ 
and on 
$S(\hat\frak h[-1]^-)$ by $t(\b_1(\!-n_1\!)\cdots \b_k(\!-n_k))\!=(-1)^k\b_1(-n_1)\cdots \b_k(-n_k)$ for $b_i\in \frak h$ and 
positive $n_i\in \frac{1}{2}+\Z.$ As before, for any $t$-stable 
subspace $W$ of $V_{\L}^T,$ we denote by $W^{\pm}$ the
eigenspaces of $t$ with eigenvalues $\pm 1.$ Then 
$(V_{\L}^T)^+$ is the tensor product of $T$ and $S(\hat\frak h[-1]^-)^-.$

The twisted vertex operator $Y(e^{\b},z)$ for $\b\in\L$ on $V_{\L}^T$ is 
defined to be
$$Y(e^{\b},z)=2^{-\<\b,\b\>}E^-_{1/2}(-\b,z)E^+_{1/2}(-\b,z)
e_{\b}z^{-\<\b,\b\>/2}$$
where
$$E^{\pm}_{1/2}(h,z)=\exp\left(\sum_{n=0}^{\infty}\frac{h(\pm (n+1/2))}{\pm
(n+1/2)}z^{\mp (n+1/2)}\right)$$ 
for $h\in \frak h,$ and $e_{\b}$ acts on $T.$ Because $\l\in2\L$ then 
 $e_{\l}$ and $e_{-\l}$
act trivially on $T,$ and we see that 
$$Y(v(\l),z)=2^{-\<\l,\l\>}E^-_{1/2}(-\l,z)E^+_{1/2}(-\l,z)z^{-\<\b,\b\>/2}
+ 2^{-\<\l,\l\>}E^-_{1/2}(\l,z)E^+_{1/2}(\l,z)z^{-\<\b,\b\>/2}$$
on $V_{\L}^T.$ 
As before we set 
$$E_{1/2}^-(\pm\l,z)E^+_{1/2}(\pm\l,z)=\sum_{n\in\Z+1/2}E^{\pm}_{1/2}(n)z^{-n}.$$
Then the contribution of $(V_{\L}^T)^+$ to $Z(v(\l),\tau)$ is
equal to 
$$q^{-1}2^{12-\<\l,\l\>}\tr (E^{+}_{1/2}(0)+E^-_{1/2}(0))q^{L(0)}|_{S(\hat\frak h[-1]^-)^-}.$$

For a formal variable $x$ we define the operator $x^N\in (\End S(\hat\frak h[-1]^-))[x]$ as before, so that $x^N(\b_1(-n_1)\cdots \b_k(-n_k))
=x^k\b_1(-n_1)\cdots \b_k(-n_k)$ for $b_i\in \frak h$ and 
positive $n_i\in \frac{1}{2}+\Z.$ Set 
$$g(q,x)=q^{3/2}\exp\left(\sum_{n=0}^{\infty}\frac{-\<\l,\l\>xq^{n+1/2}}
{(n+1/2)(1-xq^{n+1/2})}\right)\prod_{n>0}(1-xq^{n-1/2})^{-24}.$$
By a proof not essentially different to that of Lemmas \ref{l4.1} and
\ref{l4.2} we find the following:

\begin{lem} The traces 
$\tr E^{+}_{1/2}(0)q^{L(0)}x^N|_{S(\hat\frak h[-1]^-)}$ and 
$\tr E^{-}_{1/2}(0)q^{L(0)}x^N|_{S(\hat\frak h[-1]^-)}$ are
the same and equal to $g(q,x).$
\end{lem}

One can easily see that
$$\tr (E^{+}_{1/2}(0)+E^-_{1/2}(0))q^{L(0)}|_{S(\hat\frak h[-1]^-)^-}
=g(q,1)-g(q,-1).$$
Next,
\begin{eqnarray*}
& &\ \ \ \ \sum_{n=0}^{\infty}\frac{xq^{n+1/2}}{(n+1/2)(1-xq^{n+1/2})}\\
& &= \sum_{n=0}^{\infty}\frac{xq^{n+1/2}}{(n+1/2)}\sum_{i=0}^{\infty}
x^iq^{i(n+1/2)}\\
& &=\sum_{i=1}^{\infty}x^i\sum_{n=0}^{\infty}
\frac{q^{i(n+1/2)}}{(n+1/2)}\\
& & = -\sum_{i=1}^{\infty}x^i\log\left(\frac{1-q^{\frac{i}{2}}}{1+q^{\frac{i}{2}}}\right).
\end{eqnarray*}
Then the contribution of $(V_{\L}^T)^+$ to $Z(v(\l),\tau)$ is equal to
\begin{eqnarray*}
& &\ \ \ \ 2^{-\<\l,\l\>+12}q^{1/2}\prod_{n=1}^{\infty}(1-q^{n-1/2})^{-24}
\prod_{i=1}^{\infty}
\left(\frac{1-q^{\frac{i}{2}}}{1+q^{\frac{i}{2}}}\right)^{\<\l,\l\>}\\
& &\ \ \ \ -2^{-\<\l,\l\>+12}q^{1/2}\prod_{n=1}^{\infty}(1+q^{n-1/2})^{-24}
\prod_{i=1}^{\infty}
\left(\frac{1-q^{\frac{i}{2}}}{1+q^{\frac{i}{2}}}\right)^{(-1)^i\<\l,\l\>}\\
& &=2^{-\<\l,\l\>+12}q^{1/2}\prod_{n=1}^{\infty}\frac{(1-q^{n})^{24}}
{(1-q^{n/2})^{24}} \prod_{i=1}^{\infty}\frac{
(1-q^{\frac{i}{2}})^{2\<\l,\l\>}}{(1-q^i)^{\<\l,\l\>}}\\
& &\ \ \ \ -2^{-\<\l,\l\>+12}q^{1/2}\prod_{n=1}^{\infty}\frac{(1-q^{n/2})^{24}
(1-q^{2n})^{24}}
{(1-q^{n})^{48}}\prod_{i=1}^{\infty}\frac{(1-q^{i})^{5\<\l,\l\>}}
{(1-q^{2i})^{2\<\l,\l\>}(1-q^{i/2})^{2\<\l,\l\>}}\\
& &=\eta(\t)^{12}\left\{\left(\frac{\Theta_2(\tau)}{2}\right)^{\<\l,\l\>-12}
-\left(\frac{\Theta_3(\tau)}{2}\right)^{\<\l,\l\>-12}\right\}.
\end{eqnarray*}
Thus we have proved 
\begin{lem}\label{l4.10}
The contribution of $(V_{\L}^T)^+$ to $Z(v(\l),\tau)$ is equal to
$$\eta(\t)^{12}\left\{\left(\frac{\Theta_2(\tau)}{2}\right)^{\<\l,\l\>-12}
-\left(\frac{\Theta_3(\tau)}{2}\right)^{\<\l,\l\>-12}\right\}.$$
\end{lem}

Theorem \ref{thmm3} is an immediate consequence of Lemmas \ref{l4.9}
and \ref{l4.10}.

\section{A generalization of Theorem 3}
\setcounter{equation}{0}

In\ \ this\ \ section\  we generalize Theorem 3 by computing explicitly the
trace function $Z(v(\l),g,\tau)$ for certain automorphism $g$ of the
Moonshine module. As before, $\L$ is the Leech lattice.
 To describe the result we first recall some facts about 
$\Aut(V^{\natural}),$ that is to say, the Monster simple group $\M.$

The centralizer of an involution in $\M$ (of type $2B$) is a quotient
of a group $\hat C,$ partially described by the following short exact
sequence: 
$$1\to Q\to \hat C\to \Aut(\L)\to 1$$ 
where $Q\cong
2^{1+24}_+$ is an extra-special group of type $+$ and order $2^{25}.$
For more information on this and other facts we use below, see [G] or
[FLM]. The group $\hat C$ acts on both $S(\hat\frak h^-)$ and $S(\frak
h[-1]^-)$ through the natural action of $\Aut(\L),$ i.e., with kernel
$Q.$ It acts on $\C[\L]$ with kernel the center $Z(Q)$ of $Q,$ and on $T$ with
kernel a subgroup of $Z(\hat C)$ of order 
2 distinct from $Z(Q).$ Then
the quotient $C$ of $\hat C$ by the third subgroup of $Z(\hat C)$ of
order $2$ acts faithfully on $\mm.$ 

Let us fix $0\ne \l\in 2\L,$ and let $H<\hat C$ be the subgroup
defined as follows:
$$1\to Q\to H\to (\Aut\,\L)_{\l}\to 1$$ 
where $ (\Aut\,\L)_{\l}$ is the subgroup of $\Aut\, \L$ which 
fixes $\l.$ We will compute $Z(v(\l),h,\tau)$ for 
$h\in H.$ The action of $h$ on $V_{\L}$ is described by a pair 
$(\xi,a)$ where $\xi\in\frac{1}{2}\L/\L$ and $a\in  (\Aut\,\L)_{\l};$
$a$ acts in the natural manner, and
$\xi$ acts via 
$$\xi: v\otimes e^{\b}\mapsto e^{2\pi i \<\xi,\b\>}v\otimes e^{\b}.$$

We let $-a$ denote the element $ta\in \Aut\,\L,$ and define 
a modified theta-function as follows:
\begin{equation}\label{5.1}
\theta_{\xi,-a}(\t)=\sum_{\gamma\in \L, a\gamma=-\gamma}e^{2\pi i \<\xi,\gamma\>}
q^{\frac{1}{2}\<\gamma,\gamma\>}.
\end{equation}
(\ref{5.1}) is a modification of the theta-series of
the sublattice of $\L$ fixed by $-a,$ and as such is 
a modular form of weight equal to one half the dimension of the $-a$ 
fixed sublattice. Finally, let $\eta_a(\tau)$ and
$\eta_{-a}(\tau)$ by the ``usual'' eta-products associated to
$a$ and $-a$ (with regard to their action on $\L$) 
 (cf. [CN], [M]). We will establish
\begin{th}\label{t5.1} Let $0\ne \l=2\alpha,$ $\a\in \L,$ and
let $h\in H$ be associated to $(\xi,a)$ as above. Then we have 
\begin{eqnarray}\label{5.2}
& &Z(v(\l),h,\t)=e^{2\pi i\<\xi,\a\>}\frac{\theta_{\xi,a}(\tau)}{\eta_{-a}(\t)}
\left(\frac{\Theta_1(\tau)}{2}\right)^{\<\l,\l\>}\nonumber\\
& &\ \ \ \ \ \ \ +\tr_T(h)\left\{\frac{\eta_{a}(\tau)}{\eta_{a}(\tau/2)}
\left(\frac{\Theta_2(\tau)}{2}\right)^{\<\l,\l\>}
-\frac{\eta_{-a}(\tau)}{\eta_{-a}(\tau/2)}\left(\frac{\Theta_3(\tau)}{2}\right)^{\<\l,\l\>}\right\}.
\end{eqnarray}
\end{th}

Note that $\eta_{-a}(\t)$ is a form of the same weight as 
$\theta_{\xi,a}(\t)$ (loc.cit.), so that (\ref{5.2})
is indeed a form of the same weight as $Z(v(\l),\tau),$ as expected.
The proof of Theorem \ref{t5.1} is a modification of that
of Theorem \ref{thmm3}.

We\ \  begin\ \ with the appropriate modification of (\ref{4.2}), concerning
the action of $Y(v(\l),z)h$ on $V_{\L}^+.$ We have, setting
$h=h(\xi,a),$
\begin{eqnarray}
& &\ \ \  Y(v(\l),z)h(v\otimes e^{\b}+tv\otimes e^{-\b})\nonumber\\
& &=Y(v(\l),z)h(\xi,1)(a(v)\otimes e^{a(\b)}+ta(v)\otimes e^{-a(\b)})
\nonumber\\
& &=e^{2\pi i \<\xi,a(\b)\>}Y(v(\l),z)(a(v)\otimes e^{a(\b)}+ta(v)\otimes e^{-a(\b)})
\nonumber\\
& &=e^{2\pi i \<\xi,a(\b)\>}\left\{
z^{\<\l,a(\b)\>}E^-(-\l,z)E^+(-\l,z)a(v)\otimes e^{\l+a(\b)}\right.\nonumber\\
& &\ \ \ \ \ \ \ +
z^{-\<\l,a(\b)\>}E^-(-\l,z)E^+(-\l,z)ta(v)\otimes e^{\l-a(\b)}\nonumber\\
& &\ \ \ \ \ \ \ +z^{-\<\l,a(\b)\>}E^-(\l,z)E^+(\l,z)a(v)\otimes e^{-\l+a(\b)}\nonumber\\
& &\ \ \ \ \ \ \ \left.+
z^{\<\l,a(\b)\>}E^-(\l,z)E^+(\l,z)ta(v)\otimes e^{-\l-a(\b)}\right\}.\label{5.3}
\end{eqnarray}

\begin{lem}\label{l5.2} We may take $\l-a(\b)=\b$ in (\ref{5.3}).
This holds if, and only if, $\a-\b=\delta$ for some $\delta\in \L$ satisfying
$-a(\delta)=\delta.$
\end{lem}

\pf We see from (\ref{5.3}) that contributions to the trace of $o(v(\l))h$ on
$V_{\L}^+$ potentially only arise when $\l+\a(\b)=\pm \b$ or
$\l-a(\b)=\pm \b.$ If $\l=\pm(a(\b)-\b)$ then $\l$ is both a commutator 
(i.e., lies in $[a,\L]$) and a fixed-point of $a$ (by hypothesis).
This leads to the contradiction that $\l=0,$ so in fact $\l+a(\b)=-\b$ or
$\l-a(\b)=\b.$ Since $\b$ and $-\b$ are essentially interchangeable
in (\ref{5.3}), we may assume that indeed 
\begin{equation}\label{5.4}
\l-a(\b)=\b.
\end{equation}
Applying $a$ to (\ref{5.4}) yields $\l-a^2(\beta)=a(\b)=\l-\b,$ so
that $a^2(\b)=\b.$ This may be written as $(a+1)(a-1)\b=0.$ Set
\begin{equation}\label{5.5}
a(\b)-\b=2\delta.
\end{equation}
Hence $a(\delta)+\delta=0,$ that is, $2\delta$ lies in the sublattice of
$\L$ fixed by $-a.$ Moreover (\ref{5.4}) and (\ref{5.5}) yield
$\l-2\b=2\delta,$ so remembering that $\l=2\a$ we get 
\begin{equation}\label{5.6}
\a-\b=\delta.
\end{equation}

On the other hand, if (\ref{5.6}) holds, application of $a$ yields
\begin{equation}\label{5.7}
\a-a(\b)=-\delta
\end{equation}
and (\ref{5.6}), (\ref{5.7}) imply that $\l-a(\b)=\b.$
\qed

\bigskip

From the Lemma and (\ref{5.3}) we see that only expressions of the form
$$e^{2\pi i \<\xi,\a-\delta\>}
z^{-\<\l,\a-\delta\>}\left(E^-(-\l,z)E^+(-\l,z)ta(v)\otimes e^{\a-\delta}
+E^-(\l,z)E^+(\l,z)a(v)\otimes e^{-\a+\delta}\right)$$
contribute to the trace, where $\delta$ ranges over the $-a$ fixed sublattice
of $\L.$

We now follow the analysis of Section 4 which follows
(\ref{4.3}). Since $a$ fixes $\l$ then the contribution from $S(\hat A^-)$
is identical to that of (\ref{4.4}). As for $S(\hat B^-),$ the operators
$E^{\pm}(0)$ are trivial, and we need to calculate 
\begin{equation}\label{5.8}
\tr q^{L(0)}ax^N|_{S(\hat B^-)}.
\end{equation}
If $x=1$ this is precisely $\eta_a(\tau)/\eta(\t),$ by definition\footnote{
This takes into account the corresponding grade-shift.}. If $x=-1$ then $ax^N$ is just
the action of $ta,$ and (\ref{5.8}) is then 
$\frac{\eta_{-a}(\t)\eta(2\t)}{\eta(\t)}.$

Combining (\ref{4.4}) and the above, we obtain the analogue of Lemma 
\ref{l4.2}, namely:
 \begin{lem}\label{l5.3} We have for $x=\pm 1,$ 
\begin{equation}\label{5.9}
\tr E^{\pm}(0)aq^{L(0)}x^N|_{M(1)}=
\exp\left(\sum_{n>0}\frac{-\<\l,\l\>xq^n}{n(1-xq^n)}\right)\eta_{xa}(\t)^{-1}.
\end{equation}
\end{lem}

Now use this, Lemma \ref{l5.2}, and the proof of Lemma \ref{l4.3}  to see that
the contribution of $V_{\L}^+$ to $Z(v(\l),h,\t)$ is equal to
\begin{equation}\label{5.10}
\sum_{\delta\in\L, -a(\delta)=\delta} e^{2\pi i\<\xi, \alpha-\delta\>}
q^{\frac{1}{2}\< \alpha-\delta, \alpha-\delta\>}
 \exp\left(\sum_{n>0}\frac{\<\l,\l\>q^n}{n(1+q^n)}\right)\eta_{-a}(\t)^{-1}.
\end{equation}

Note that $\<\a,\delta\>=0.$ Then (\ref{5.10}) is equal to 
$$e^{2\pi i\<\xi, \alpha\>}\theta_{\xi,-a}(\t)
\left(\frac{\Theta_1(\t)}{2}\right)^{\<\l,\l\>}\eta_{-a}(\t)^{-1}$$
which is the first summand of (\ref{5.2}). 

The other two summands of (\ref{5.2}) arise from the contribution of 
$(V_{\L}^T)^+$ to the trace. The proofs are as before, and are easier
than the part just completed as there is no
theta-function to deal with. We leave details to the reader. This
completes our discussion of Theorem \ref{t5.1}.

\section{The invariance group of $Z(v,g,\t)$}
\setcounter{equation}{0}

We will determine the subgroup of $\Gamma=SL(2,\Z)$ which leaves $Z(v,g,\tau)$
invariant. More precisely,
if $v$ is homogeneous of weight k with respect to $L[0],$ so that $Z(v,g,\t)$
is modular of weight $k$ by [DLM], we
will describe in Theorem \ref{t6.1} below the action of $\Gamma_0(n)$
on $Z(v,g,\t),$ where $n$ is the order of $g.$

The case where $v=\bc$ is the vacuum (and $k=0$) is covered by results
in [CN] and [B2]. Precisely, one knows that there is a character $\e_g$ of
$\Gamma_0(n)$ such that 
\begin{equation}\label{6.1}
Z|\gamma(\bc,g,\tau):=Z(\bc,g,\gamma\t)=\e_g(\gamma)Z(\bc,g,\t)
\end{equation}
for $\gamma\in\Gamma_0(n).$ Moreover $\ker\e_g=\Gamma_0(N)$
where $N=nh,$ and $h$ divides $gcd(n,24).$

To describe our generalization of this result, we need to recall some
further results. Let $A_{\M}(\<g\>)=N_{\M}(\<g\>)/ C_{\M}(\<g\>)$ be
the {\em automizer} of $\<g\>$ in the Monster $\M.$ Then
$A_{\M}(\<g\>)$ is the group of automorphisms of $\<g\>$ induced by
conjugation in $\M.$ As such, $A_{\M}(\<g\>)$ has a canonical
embedding
\begin{equation}\label{6.2}
i_g: A_{\M}(\<g\>)\to U_n
\end{equation}
in which $U_n$ is the group of units of $\Z/n\Z,$ and $t\in
N_{\M}(\<g\>)$ satisfying $tgt^{-1}=g^d$ maps to $d$ under $i_g.$
From the character table of $\M$ [Cal], we see that the following is true:
$[U_n:\im i_g]\leq 2,$ with equality if, and only if, $g$ is 
{\em not} conjugate to $g^{-1}$ in $\M.$ In this case, $U_n=\im i_g\times
\{\pm 1\}.$

Since $\Gamma_0(n)/\Gamma_1(n)$ is naturally isomorphic to 
$U_n,$ we may define a subgroup $\Gamma_g(n)$ of $\Gamma_0(n)$
via the following diagram (rows being short exact)
\begin{equation}\label{6.3}
\begin{array}{ccccccccc}
1 &\to& \Gamma_1(n)&\to&\Gamma_0(n)&\to &U_n& \to& 1\\
& &\|&&\uparrow&&\uparrow&&\\
 1 &\to& \Gamma_1(n)&\to&\Gamma_g(n)&\to &i_g(A_{\M}(\<g\>))& \to& 1
\end{array}
\end{equation}
In (\ref{6.3}), $\gamma=\left(\begin{array}{cc} a & b\\
c & d
\end{array}
\right)\in \Gamma_0(n)$ maps to $d\in U_n.$
From what we have said, we have
$[\Gamma_0(n):\Gamma_g(n)]\leq 2,$ and $\Gamma_0(n)=\Gamma_g(n)\times \{\pm 1\}$
if we have equality.

Let $\chi$ range over the irreducible, complex characters of the normalizer 
$N_{\M}(\<g\>)$ of $\<g\>$ in $\M.$ We will be particularly interested 
in those $\chi$ satisfying $C_{\M}(\<g\>)\subset \ker \chi.$ Such
$\chi$ are 1-dimensional, and induce characters
\begin{equation}\label{6.4}
\chi:A_{\M}(\<g\>)\to \C^*.
\end{equation}
Using the lower row of (\ref{6.3}), we can pull-back $\chi$ to a 
character of $\Gamma_g(n),$ also denoted by $\chi.$ If 
$[\Gamma_0(n):\Gamma_g(n)]= 2$ then  $\Gamma_0(n)=\Gamma_g(n)\times 
\<-I\>$ (where $I$ is the $2\times 2$ identity matrix) 
and we then define a character $\chi_k$ $(k\in\Z)$ of $\Gamma_0(n)$
so that its restriction to $\Gamma_g(n)$ is the earlier
$\chi,$ and its value on $-I$ is $(-1)^k.$ So in all cases we have defined 
 characters
$\chi_k$ of $\Gamma_0(n),$ with 
the convention that $\chi_k=\chi$ if $\Gamma_g(n)=\Gamma_0(n).$

We decompose $\mm$ into homogeneous subspaces $\mm_{[k]}$ with 
respect to the $L[0]$-operator. This
commutes with the action of the Monster $\M,$ and we let $\mm_{[k],\chi}$
be the $\chi$-isotypic subspaces of $\mm_{[k]}$ considered as a 
$N_{\M}(\<g\>)$-module. We can now state our result:
\begin{th}\label{t6.1}
Fix $g\in\M,$ and let the notation be as above. Suppose that
$v\in \mm_{[k],\chi}$ for some simple character $\chi$ of
$N_{\M}(\<g\>).$ Then the following hold:

(a) If $C_{\M}(\<g\>)\not\subset \ker \chi$ then $Z(v,g,\tau)=0.$

(b) If  $C_{\M}(\<g\>)\subset \ker \chi$ then
\begin{equation}\label{6.5}
Z|_k\gamma(v,g,\tau)=\e_g(\gamma)\overline{\chi_k(\gamma)}Z(v,g,\tau)
\end{equation}
for $\gamma\in \Gamma_0(n).$
\end{th}

\pf We first prove (a). Since $\chi$ is a simple character of
$N_{\M}(\<g\>)$ and $C_{\M}(\<g\>)$ is normal in $N_{\M}(\<g\>),$
the assumption $C_{\M}(\<g\>)\not\subset \ker \chi$ means that
$C_{\M}(\<g\>)$  does not leave $v$ invariant if $0\ne v\in \mm_{[k],\chi}.$
Then $v$ can be written as a linear combination $v=\sum_{i}v_i$
with each $v_i\in \mm_{[k],\chi}$ and $t_iv_i=\l_iv_i$
for each $i,$ some $t_i\in C_{\M}(\<g\>),$ and $1\ne \l_i\in \C^*.$

We may thus assume that $v=v_i,$ with $tv=\l v$ for some $t\in C_{\M}(\<g\>)$
and some $1\ne \l\in \C^*.$ But
then 
\begin{eqnarray*}
& &Z(v,g,\t)=q^{-1}\sum_{n}(\tr|_{\mm_n}o(v)g)q^n\\
& &\ \ \ \ =q^{-1}\sum_{n}(\tr|_{\mm_n}to(v)gt^{-1})q^n\\
& &\ \ \ \ =q^{-1}\sum_{n}(\tr|_{\mm_n}o(tv)g)q^n\\
& &\ \ \ \ =Z(tv,g,\t)\\
& &\ \ \ \ =\l Z(v,g,\tau).
\end{eqnarray*}
Since $\l\ne 1,$ we get $Z(v,g,\t)=0,$ as required.

To prove (b) we need some results from [DLM], which we assume the 
reader is familiar with. In particular, since $g$ has order $n$ then
a matrix $\gamma=\left(\begin{array}{cc}
a & b\\
c & d
\end{array}
\right)\in\Gamma_0(n)$ maps the $(1,g)$ conformal block to the $(1,g^d)$
conformal block. Since the trace functions $Z(v,g,\tau),$
$Z(v,g^d,\t)$ span these conformal blocks, there is a scalar
$\eta_g(\gamma),$ independent of $v,$ such that
\begin{equation}\label{6.6}
Z|\gamma(v,g,\tau)=\eta_g(\gamma)Z(v,g^d,\t).
\end{equation}
Here, if $v\in  \mm_{[k]}$ then 
\begin{equation}\label{6.7}
Z|\gamma(v,g,\tau)=(c\t+d)^{-k}Z(v,g,\gamma\t).
\end{equation}
Taking $v=\bc, k=0$ in (\ref{6.6})-(\ref{6.7}) and comparing with (\ref{6.1})
then yields $\eta_g(\gamma)=\e_g(\gamma),$ that is
\begin{equation}\label{6.8}
Z|\gamma(v,g,\tau)=\e_g(\gamma)Z(v,g^d,\t).
\end{equation}
Suppose that $d\in i_g(A_{\M}(\<g\>),$ that is $\gamma\in \Gamma_g(n).$
Then $g^d=tgt^{-1}$ for some $t\in N_{\M}(\<g\>),$ and we calculate as
before:
$$Z(v,g^d,\t)=Z(v,tgt^{-1},\t)=Z(t^{-1}v,g,\t)=\chi(t^{-1})Z(v,g,\t).$$
Then (\ref{6.8}) reads 
\begin{equation}\label{6.9}
  Z|\gamma(v,g,\tau)=\e_g(\gamma)\chi(t^{-1})Z(v,g,\t).
\end{equation}
By our conventions, $\chi(t^{-1})=\overline{\chi(\gamma)},$ so 
(\ref{6.9}) is what we require.

Now assume that $\gamma\not \in \Gamma_g(n).$ From our earlier
remarks, it suffices to take $\gamma=-I.$ In this case $\gamma\in \Gamma_0(N),$ so $\e_g(\gamma)=1,$ and
(\ref{6.7}) reads
$$Z|_k\gamma(v,g,\t)=(-1)^kZ(v,g,\t),$$
which is what (\ref{6.5}) says in this case. The proof of theorem
is now complete. \qed

\begin{rem} {\rm By Theorem 2 of [DM], each $\chi$ occurs in $\mm,$
that is, given $\chi$ as above, there is $k$ such that $V_{[k],\chi}\ne 0.$}
\end{rem}

\end{document}